\documentclass[preprint,1p]{elsarticle} 

\usepackage{graphicx,bbm,url,float,booktabs} 
\usepackage[colorlinks=true,linkcolor=black, citecolor=black]{hyperref}
\usepackage{algorithmic,algorithm,amsfonts} 
\usepackage[fleqn]{amsmath}
\usepackage{algorithmic}
\usepackage{algorithm}

\newcommand{\ignore}[1]{ }
\newcommand{\R}[0]{R}

\begin{document}
\begin{frontmatter} 
\title{A Space-Indexed Formulation of Packing Boxes into a Larger Box}
\tnotetext[label1]{All accompanying materials are available at \url{http://discretisation.sourceforge.net}.}

\author{Sam D. Allen}
\author{Edmund K. Burke}
\author{Jakub Mare\v{c}ek \corref{cor1}}
\ead{jakub@marecek.cz}
\cortext[cor1]{The corresponding author is Jakub Mare\v{c}ek.}
\address{The research has been performed, while all authors were at the University of Nottingham, School of Computer Science, Jubilee Campus, Nottingham, NG8 1BB, UK
} 



\date{Received: \today
}

\begin{abstract}
Integer programming solvers fail to decide whether 12 unit cubes can be packed into a 1x1x11 box within an hour using the natural relaxation of Chen/Padberg. We present an alternative relaxation of the problem of packing boxes into a larger box, which makes it possible to solve much larger instances.
\ignore{
The problem of packing many small boxes into a single larger box underlies a number of cutting, packing, scheduling, and transportation applications. 
There are a number of heuristic solvers, but the progress in exact solvers, in general, and integer programming solvers, in particular, has been limited. 
Padberg [Math. Methods Oper. Res., 52(1):1--21, 2000] estimated his extension of the integer linear programming formulation of Chen et al.\ could 
cope with ``about twenty boxes''. 
Our computational experiments confirm that the seemingly trivial decision whether twelve unit-cubes can be packed into a box with the unit-base and 
height eleven (``3D Pigeon Hole''), cannot be made within an hour by modern integer programming solvers using this formulation.  

We present a new, ``space-indexed'', linear programming relaxation, which often provides lower bounds within 1 percent of optimality, 
and makes it possible to solve instances of 3D Pigeon Hole problem with ten million boxes within an hour.
Results of extensive computational tests of both formulations are reported.  
}
\end{abstract}
\end{frontmatter}

\section{Introduction}

Multi-dimensional packing problems have a number of applications in cutting, packing, scheduling, and transportation.
Problems in dimension three with rotations around combinations of axes in multiples of 90 degrees
are of particular interest in many natural applications.
Let us fix the order of six such allowable rotations in dimension three arbitrarily and define:

\vskip 6pt
\noindent
The {\sc Container Loading Problem (CLP)}: Given dimensions of a large box (``container'') $x, y, z > 0$ and dimensions of
$n$ small boxes $D \in {\mathbbm{R}}^{n \times 3}$ with associated values $w \in {\mathbbm{R}}^{n}$, 
and specification of the allowed rotations $r = \{0, 1\}^{n \times 6}$, find the 
greatest $k \in \mathbbm{R}$ such that there is a packing of small boxes $I \subseteq \{1, 2, \ldots, n\}$ 
into the container 
with value $k = \sum_{i \in I} w_i$. 
The packed small boxes $I$ may be rotated in any of the allowed ways,
must not overlap, and no vertex can be outside of the container.

\vskip 6pt
\noindent
The {\sc Van Loading Problem (VLP)}: Given dimensions of a large box (``van'') $x, y, z > 0$, maximum mass $p \ge 0$ it can hold (``payload''), 
dimensions of $n$ small boxes $D \in {\mathbbm{R}}^{n \times 3}$ with associated values $w \in {\mathbbm{R}}^{n}$, 
mass $m \in {\mathbbm{R}}^{n}$, 
and specification of the allowed rotations $r = \{0, 1\}^{n \times 6}$, find the greatest $k \in \mathbbm{R}$ 
such that there is a packing of small boxes $I \subseteq \{1, 2, \ldots, n\}$ into the container 
with value $k = \sum_{i \in I} w_i$ and mass $\sum_{i \in I} m_i \le p$.
The packed small boxes $I$ may be rotated in any of the allowed ways,
must not overlap, and no vertex can be outside of the container.
\vskip 6pt

Both problems are NP-Hard even to approximate \cite{Chlebik09}.  
For example in branch-and-price-and-cut solvers for realistic models of transportation considering both weight and volume of vehicle load \cite{Gendreau06,LuebbeckePC},
however, one needs to solve the van loading problem to optimality as the pricing subproblem.
In scheduling, the problem corresponds to scheduling non-malleable parallel jobs, 
requiring a certain number of clock cycles, a certain number of processors, and a certain amount of another discrete resource.
Although there are a number of excellent heuristic solvers,
the progress in exact solvers for the Container Loading Problem has been limited, so far.  
Chen et al.\ \cite{Chen1995} proposed an integer linear programming (ILP) formulation, 
which allowed for up to six boxes to be packed into a smaller box.
Fasano \cite{Fasano99} and Padberg \cite{Padberg00} improved the formulation and suggested it could cope with ``about twenty boxes''. 
We confirm these limits using ILOG CPLEX, Gurobi, and SCIP in Section~\ref{sec:computational}.
This observation motivates the rest of the paper, where we present a new space-indexed formulation providing 
strong linear programming relaxations. 
Throughout, we use the notation as described in Table~\ref{notationtable}. 


\begin{table}[t]
	\centering
		\begin{tabular}{cl}
		\toprule \\
		Symbol & Meaning \\
		\midrule
			$n$ & The number of boxes. \\
			$H$ & A fixed axis, in the set $\{X,Y,Z\}$. \\
			$\alpha$ & An axis of a box, in the set $\{1,2,3\}$. \\
			$L_{\alpha i}$ & The length of axis $\alpha$ of box $i$. \\
			$l_{\alpha i}$ & The length of axis $\alpha$ of box $i$ halved. \\
			$D_H$ & The length of axis $H$ of the container. \\
			$w_i$ & The volume of box $i$ in the CLP. \\ 
     \bottomrule
		\end{tabular}
	\caption{Notation used.}
	\label{notationtable}
\end{table}

\section{Related Work}

Chen et al.\ \cite{Chen1995} introduced an integer linear programming formulation using the relative placement indicator:

\begin{align}
\lambda_{ij}^H = \begin{cases}
   \; 1 & \text{if box } i \text{ precedes box } j \text{ along axis } H \\
   \; 0 & \text{otherwise} \notag 
 \end{cases}
\end{align}
\begin{align}
\delta_{\alpha i}^H = \begin{cases}
   \; 1 & \text{if box } i \text{ is rotated so that axis } \alpha \text{ is parallel to fixed } H \\
   \; 0 & \text{otherwise} \notag 
 \end{cases} 
\end{align}

Using implicit quantification, it reads: 

\begin{flalign}
\max &\sum_{i=1}^n{\sum_H{w_i\delta_{1i}^H}}\\
\textrm{s.t. } &\sum_H{\delta_{2i}^H} = \sum_H{\delta_{1i}^H} = \sum_\alpha \delta_{\alpha i}^H \\ 
&L_{1j(i)}\lambda^H_{j(i)i} + \sum_\alpha l_{\alpha i}\delta_{\alpha i}^H \leq x_i^H \leq \sum_\alpha{(D_H - l_{\alpha i})\delta_{\alpha i}^H} - L_{1j(i)}\lambda_{ij(i)}^H \\ 
&D_H\lambda_{ji}^H + \sum_\alpha{l_{\alpha i}\delta_{\alpha i}^H} - \sum_\alpha{(D_H - l_{\alpha j})\delta_{\alpha j}^H} \leq x_i^H - x_j^H \\ 
&\leq \sum_\alpha{(D_H - l_{\alpha i})\delta_{\alpha i}^H} - \sum_\alpha{l_{\alpha j}\delta_{\alpha j}^H} - D_H\lambda_{ij}^H\\ 
&\sum_H{(\lambda_{ij}^H + \lambda_{ji}^H)} \leq \sum_H{\delta_{1i}^H}, \sum_H{(\lambda_{ij}^H + \lambda_{ji}^H)} \leq \sum_H{\delta_{1j}^H} \\ 
&\sum_H{\delta_{1i}^H} + \sum_H{\delta_{1j}^H} \leq 1 + \sum_H{(\lambda_{ij}^H + \lambda_{ji}^H)}\\ 
&\sum_{i=1}^n{\sum_H{\left(\prod_\alpha {L_{\alpha i}}\right)\delta_{1i}^H \leq \prod_H{D_H}}}\\ 
&\delta_{\alpha i}^H \in \{0,1\}, \lambda_{ij}^H \in \{0,1\} \\ 
& 
L_{1i}\leq L_{2i} \leq L_{3i}, j(i) \text{ such that } L_{1j(i)} = \max\{L_{1j}\} \text{ for } 1 \leq i \neq j \leq n.
\label{clpmodel}
\end{flalign}

Padberg \cite{Padberg00} has studied its properties. In particular, he identified
the subsets of constraints with the integer property. 
Despite the interesting theoretical properties, 
modern integer programming solvers fail to solve instances 
larger than 10--20 boxes using this formulation,
as evidenced by Table~\ref{tab:pigeon-results}.
This is far from satisfactory.

\subsection{Extensions to the Formulation of Chen/Padberg}
\label{extensions} 

Arguably, the formulation could be strengthened by addition of further valid constraints.
We have identified three classes of such constraints.
Compound constraints stop combinations of certain boxes being placed next to each other if they would violate the domain constraint.
For example:
\begin{eqnarray}
\delta_{ki}^H + \delta_{mj}^H + \lambda_{ij}^H \leq 2 & \text{ if } l_{ki} + l_{mj} > D_{H} & \forall k,m\in\{1,2,3\}, H\in\{X,Y,Z\} \notag\\ \mbox{ and } 1\leq i \neq j \leq n
\label{compoundconstraints}
\end{eqnarray}
The example for 2 boxes can be easily extended to 3 or more boxes ``in a row.''
One can also attempt to break symmetries in the problem.
If boxes $i$ and $j$ (where $i < j$) are identical (i.e. $L_{\alpha i}=L_{\alpha j}~\forall \alpha$):
\begin{eqnarray}
\lambda_{ij}^H = 0 \mbox{ and } \sum_H{\delta_{1i}^H} \geq \sum_H{\delta_{1j}^H}
\label{identicalboxsymmetrybreaking2}
\end{eqnarray}
Furthermore, if a box has two or more sides of the same length then we can limit the rotations.
We define function $f$, which provides a canonical mapping of $\{X,Y,Z\}$ to $\{1,2,3\}$:
\begin{eqnarray}
\sum_H{(f(H) \cdot \delta_{ki}^H)} \geq \sum_H{(f(H) \cdot \delta_{mi}^H)} & \text{if } l_{ki} = l_{mi} & \forall 1\leq i \leq n \text{ and } 1 \leq k < m \leq 3
\label{rotationalsymmetrybreaking}
\end{eqnarray}
Nevertheless, whilst the addition of these constraints improves the performance somewhat, 
the formulation remains impractical.

\section{A Space-Indexed Formulation}
\label{sec:formulation}

We hence propose a new formulation, where the large box is discretised into units of space, whose dimensions are given by
the largest common denominator of the respective dimensions of the small boxes.
The small boxes are grouped by their dimensions into types $t \in \{1,2,\ldots,n\}$. 
$A_t$ is the number of boxes of type $t$ available.
The formulation uses the space-indexed binary variable:

\begin{align}
\mu_{x,y,z}^{t}  = \begin{cases}
   \; 1 & \text{if a box of type } t \text{ is placed such that its lower-back-left} \\
        & \text{ vertex is at coordinates } x, y, z \\
   \; 0 & \text{otherwise} 
 \end{cases}
\label{eq:muvars}
\end{align}

Without allowing for rotations, the formulations reads:

\begin{align}
\max & \sum_{x,y,z,t}{\mu_{xyz}^t w_t} \label{objective-discretised} \\
\text{s.t. } \sum_t{\mu_{xyz}^t} & \leq 1 \quad \forall x,y,z \label{muconstraint1} \\
\mu_{xyz}^t & = 0 \quad \forall x,y,z,t \text{ where } x+L_{1t} > D_X \text{ or } \notag \\
            & \quad \quad \text{ or } y+L_{2t} > D_Y \text{ or } z+L_{3t} > D_Z \label{containedconstraint}\\
\sum_{x,y,z}{\mu_{xyz}^t} & \leq A_t \quad \forall t \label{muconstraint2} \\
\sum_{x',y',z',t' \in f(D, n, L, x, y, z, t)} \mu_{x'y'z'}^{t'} &\leq 1 \quad \forall x,y,z,t \label{no-overlap}\\
\mu_{xyz}^t & \in \{ 0, 1\} \quad \forall x,y,z,t
\end{align}
where one may use Algorithm~\ref{algo:genNonOverlap} in the Appendix or similar\footnote{
For one particularly efficient version of Algorithm~\ref{algo:genNonOverlap}, 
see class {\tt DiscretisedModelSolver} and its method {\tt addPositioningConstraints} in file 
{\tt DiscretisedModelSolver.cpp}
available at \url{http://discretisation.sourceforge.net/spaceindexed.zip} (September 30th, 2014).
This code is distributed under GNU Lesser General Public License.
}
 to generate the index set $f$ in Constraint~\ref{no-overlap}.

The constraints are very natural: No region in space may be occupied by more than one box type (\ref{muconstraint1}),
boxes must be fully contained within the container (\ref{containedconstraint}),
there may not be more than $A_t$ boxes of type $t$ (\ref{muconstraint2}),
and boxes cannot overlap (\ref{no-overlap}). 
There is one non-overlapping constraint (\ref{no-overlap}) for each discretised unit of space and type of box.
In order to support rotations, new box types need to be generated for each allowed rotation and 
linked via set packing constraints similar to Constraint~\ref{muconstraint2}.
In order to extend the formulation to the van loading problem, 
it suffices to add the payload capacity constraint: 

\begin{align}
	\sum_{x,y,z,t}{\mu_{xyz}^t m_t} \leq p
	\label{payloadconstraint-discretised}
\end{align}

\subsection{Adaptive Discretisation}

In order to reduce the number of regions of space, and thus the number of variables 
in the formulation, a sensible space-discretisation method should be employed. 
Mare{\v c}ek \cite{Marecek2011} has shown the problem of finding the best possible 
discretisation is $\Delta_2^p$-Complete, where $\Delta_2^p$ is the class of problems 
that can be solved in polynomial time using oracles from NP.
There are, however, practical heuristics.

The greatest common divisor (GCD) reduction can be applied on a per-axis basis, 
finding the greatest common divisor between the length of the container for an axis 
and all the valid lengths of boxes that can be aligned along that axis and scaling 
by the inverse of the GCD. 
This is trivial to do and is useful when all lengths are multiples of each other, 
which may be common in certain classes of instances.
In other instances, this may not reduce the number of variables at all. 
In such cases, a non-linear space-discretisation approach can be applied. 
To do this we use the same values as before on a per-axis basis, i.e. the lengths 
of any box sides that can be aligned along the axis. 
We then use dynamic programming to generate all valid locations for a box to be placed. 
For example, given an axis of length 10 and box lengths of 3, 4 and 6, we can place 
boxes at positions 0, 3, 4, 6, and 7.
(Positions 8, 9 and 10 are of course also possible, 
but no length is small enough to still lie within the container if placed at these points.) 
This has reduced the number of regions along that axis from 10 to 5. 
An improvement on this scale may not be particularly common in practice, though the approach can help where the GCD is 1. It is obvious that this approach can be no worse than the GCD method at discretising the container. 
This also adds some implicit symmetry breaking into the model. 

\section{Computational Experience}
\label{sec:computational}

We have tested the formulations above on three sets of instances:
\begin{itemize}
\item 3D Pigeon Hole Problem instances, Pigeon-$n$, where $n + 1$ unit cubes are to be packed
into a container of dimensions $(1 + \epsilon) \times (1 + \epsilon) \times n$.\footnote{For $n=10, 11, 12, 13, 19$, these instances are included in the MIPLIB 2010, cf. \url{http://miplib2010.zib.de/}.}
\item SA and SAX datasets, which are used to test the dependence of solvers' performance 
on parameters of the instances, notably the number of boxes, heterogeneity of the boxes, and physical dimensions of the container.
There is 1 pseudo-randomly generated instance for every combination of container sizes ranging from 5--100 in steps of 5 
units cubed and the number of boxes to pack ranging from 5--100 in steps of 5. 
The SA datasets are perfectly packable, i.e. are guaranteed to be possible to load the container 
with 100\% utilisation with all boxes packed. 
The SAX are similar but have no such guarantees; the summed volume of the boxes is greater than that of the container. 
\item Van-loading instances, generated so as to be similar to those used in pricing subproblems of branch-and-cut-and-price approaches to vehicle 
routing with both load weight and volume considerations \cite{Gendreau06,LuebbeckePC}. 
The van-loading and Van-loading-X instances use containers of 10x10x30 and 10x10x50 units, respectively, 
representing the approximate ratios of a small commercial van and a larger freight truck.
2, 4 or 6 out of 6 pseudo-randomly chosen orientations are allowed. 
The load weights and knapsack values are generated pseudo-randomly and independently of the volume of the small boxes. 
10 instances were generated for each class of problem.
\end{itemize}

All the tests were performed on a 64-bit computer running Linux, 
 which was equipped with 2 quad-core processors (Intel Xeon E5472) and 16~GB memory. 
The solvers tested were IBM ILOG CPLEX 12.2.0, Gurobi Solver 4.0, and SCIP 2.0.1 \cite{Achterberg2007} 
with CLP  as the linear programming solver. The instances are available at \url{http://discretisation.sourceforge.net}.

For the 3D Pigeon Hole Problem, results obtained within one hour using the three solvers and the Chen/Padberg formulation are shown in 
Table~\ref{tab:pigeon-results}, while Table~\ref{tab:pigeontwo} compares the results of Gurobi Solver on both formulations. 
None of the solvers managed to prove optimality of the incumbent solution for twelve unit-cubes within an hour using the Chen/Padberg formulation, 
although the instance of linear programming (after pre-solve) had only 616 rows and 253 columns and 4268 non-zeros. 
Due to their surprising difficulty, Chen/Padberg formulation of Pigeon-$n$ for $n \in \{ 10, 11, 12, 13, 19 \}$ have been contributed to and accepted into MIPLIB 2010 \cite{MIPLIB2010}, the benchmark.
In contrast, the space-indexed formulation allowed Gurobi Solver to solve Pigeon-10000000 within an hour, 
where the instance of linear programming had 10000002 rows, 10000000 columns, and 30000000 non-zeros.

\begin{figure*}[h]
	\centering
	\includegraphics[scale=0.4]{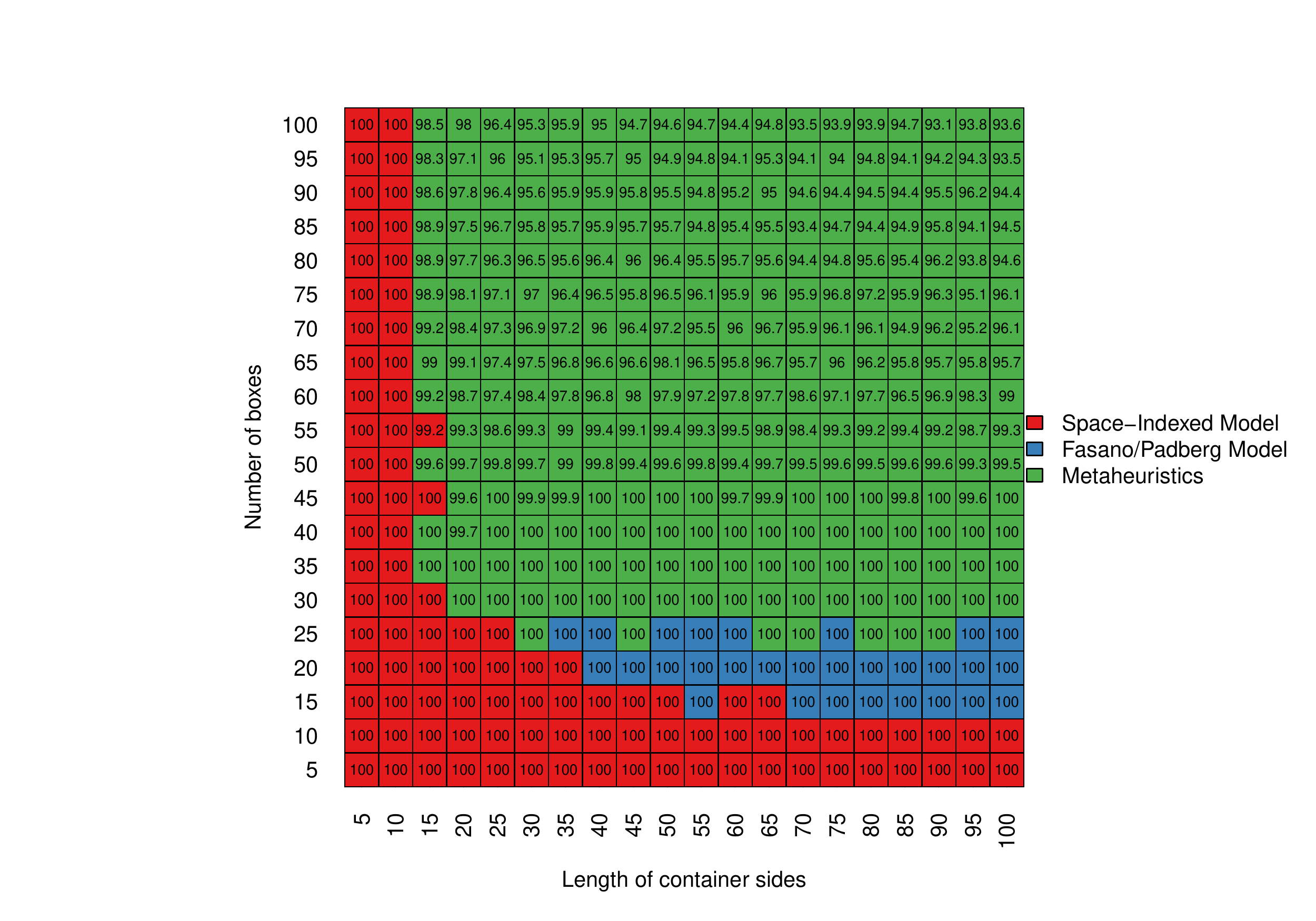}
	\caption{The best solutions obtained within an hour per solver per instance from the SA dataset. 
  Each square represents 1 pseudo-randomly generated instance.
  The number is the volume utilisation (in percent) with the tight upper bound of 100.}
	\label{contourplot}
\end{figure*}

\begin{figure*}[h]
	\centering
	\includegraphics[scale=0.4]{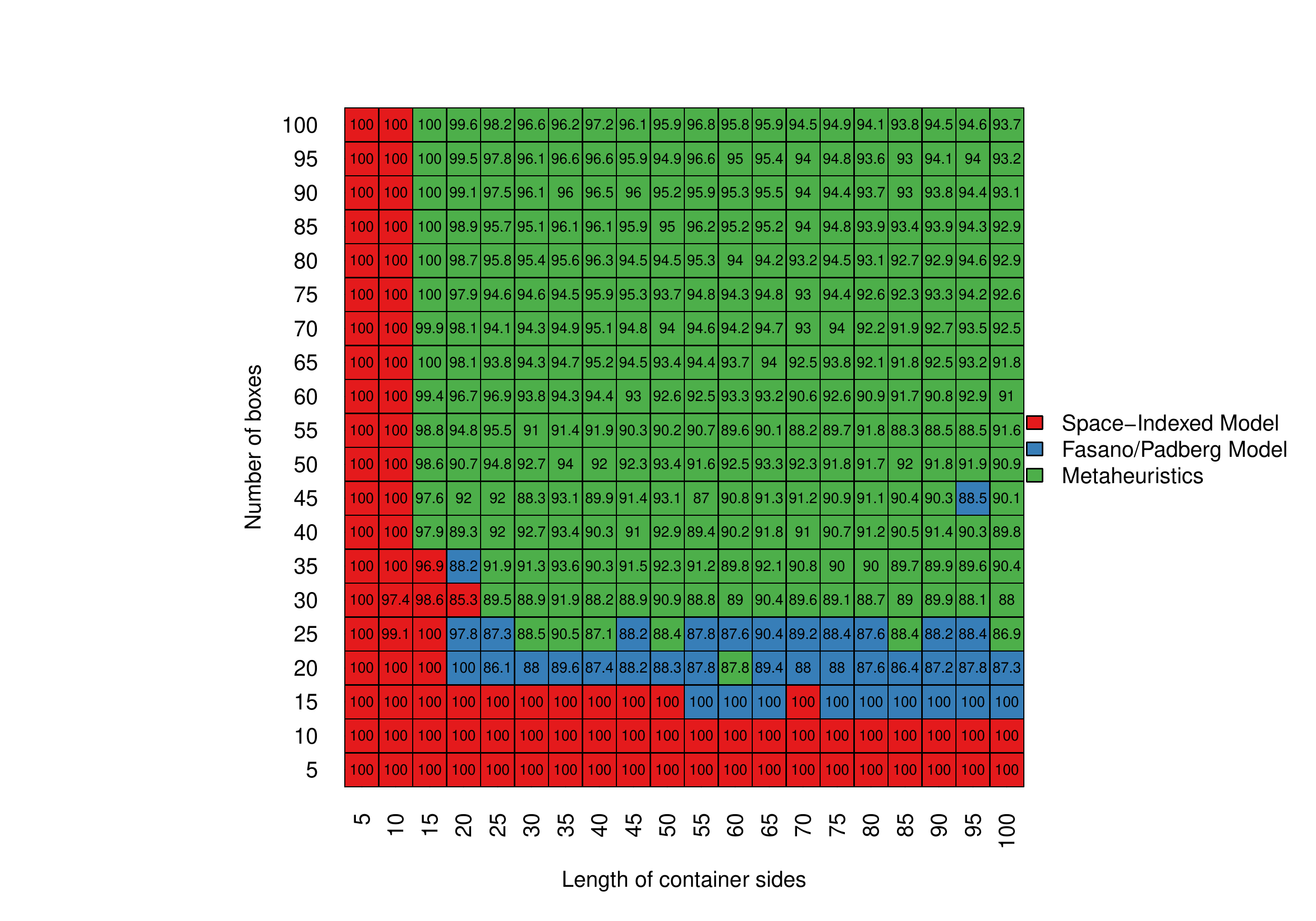}
	\caption{The best solutions obtained within an hour per solver per instance from the SAX dataset.
  Each square represents 1 pseudo-randomly generated instance.
  The number displayed is $100 (1 - s/b)$ for solution with value $s$ and upper bound $b$.
  }
	\label{contourplotX}
\end{figure*}

For the SA and SAX datasets, Figures \ref{contourplot}--\ref{contourplotX} summarise solutions obtained within an hour using the 
Chen/Padberg formulation, the space-indexed formulation, and a metaheuristic approach described by Allen et al.\ \cite{Allen11,Allen2011}.
In the case of the SAX dataset, the volume utilisation is given with respect to the tightest upper bound found by any of the solvers.
Finally, for the van-loading instances, the results of Gurobi Solver with the space-indexed model after one hour can be seen in Table \ref{tab:vanloading}. 
As there were 10 instances tested for each class of problem,
the mean number of boxes, number of unique box types and mean gap between solution value and bound are also given.

Overall, the space-indexed relaxation provides a particularly strong upper bound.
The mean integrality gap, or the ratio of the difference between root linear programming relaxation value and optimum to optimum
has been 10.49 \% and 0.37 \% for the Chen/Padberg and the space-indexed formulation, respectively, on the SA and SAX instances solved to optimality 
within the time limit of one hour. 
The mean gap, or the ratio of the difference between linear programming relaxation value and value of the best solution found to the bound, 
left after an hour on Van-loading instances was 1.5\%.
This is encouraging, albeit perhaps not particularly surprising, as similar discretised relaxations proved to be very strong 
in scheduling problems corresponding to one-dimensional packing \cite{Sousa1992,MR1308272,MR2324961} and can be shown to be 
asymptotically optimal for various geometric problems both in dimensions two \cite{MR623065} and higher \cite{MR777998} dimensions. 

\begin{table}[h]
  \centering
  \caption{The performance of various solvers on 3D Pigeon Hole Problem instances encoded in the Chen/Padberg formulation. 
  ``--'' denotes that optimality of the incumbent solution has not been proven within an hour.}
    \begin{tabular}{lrrr}
    \addlinespace
    \toprule
          &       & Time (s) &  \\
    \midrule
          & Gurobi 4.0 & CPLEX 12.2.0 & SCIP 2.0.1 + CLP \\
    Pigeon-01 & \multicolumn{1}{c}{$<$ 1} & \multicolumn{1}{c}{$<$ 1} & \multicolumn{1}{c}{$<$ 1} \\
    Pigeon-02 & \multicolumn{1}{c}{$<$ 1} & \multicolumn{1}{c}{$<$ 1} & \multicolumn{1}{c}{$<$ 1} \\
    Pigeon-03 & \multicolumn{1}{c}{$<$ 1} & \multicolumn{1}{c}{$<$ 1} & \multicolumn{1}{c}{$<$ 1} \\
    Pigeon-04 & \multicolumn{1}{c}{$<$ 1} & \multicolumn{1}{c}{$<$ 1} & \multicolumn{1}{c}{$<$ 1} \\
    Pigeon-05 & \multicolumn{1}{c}{$<$ 1} & \multicolumn{1}{c}{$<$ 1} & \multicolumn{1}{c}{3.3} \\
    Pigeon-06 & \multicolumn{1}{c}{$<$ 1} & \multicolumn{1}{c}{$<$ 1} & \multicolumn{1}{c}{37.9} \\
    Pigeon-07 & \multicolumn{1}{c}{1.5} & \multicolumn{1}{c}{3.6} & \multicolumn{1}{c}{779.3} \\
    Pigeon-08 & \multicolumn{1}{c}{7.4} & \multicolumn{1}{c}{25.6} & \multicolumn{1}{c}{--} \\
    Pigeon-09 & \multicolumn{1}{c}{88.6} & \multicolumn{1}{c}{398.4} & \multicolumn{1}{c}{--} \\
    Pigeon-10 & \multicolumn{1}{c}{1381.4} & \multicolumn{1}{c}{--} & \multicolumn{1}{c}{--} \\
    Pigeon-11 & \multicolumn{1}{c}{--} & \multicolumn{1}{c}{--} & \multicolumn{1}{c}{--} \\
    Pigeon-12 & \multicolumn{1}{c}{--} & \multicolumn{1}{c}{--} & \multicolumn{1}{c}{--} \\
    \bottomrule
    \end{tabular}%

  \label{tab:pigeon-results}
\end{table}

\begin{table}[h]
  \centering
  \caption{The performance of Gurobi 4.0 on 3D Pigeon Hole Problem instances encoded in the Chen/Padberg 
  and the space-indexed formulations. ``--'' denotes no integer solution has been found. }
    \begin{tabular}{lrr}
    \addlinespace
    \toprule
          & \multicolumn{2}{c}{Time (s)} \\
    
    \midrule
          & Chen/Padberg & Space-indexed \\
    Pigeon-01 & \multicolumn{1}{c}{$<$ 1} & \multicolumn{1}{c}{$<$ 1} \\
    Pigeon-02 & \multicolumn{1}{c}{$<$ 1} & \multicolumn{1}{c}{$<$ 1} \\
    Pigeon-03 & \multicolumn{1}{c}{$<$ 1} & \multicolumn{1}{c}{$<$ 1} \\
    Pigeon-04 & \multicolumn{1}{c}{$<$ 1} & \multicolumn{1}{c}{$<$ 1} \\
    Pigeon-05 & \multicolumn{1}{c}{$<$ 1} & \multicolumn{1}{c}{$<$ 1} \\
    Pigeon-06 & \multicolumn{1}{c}{$<$ 1} & \multicolumn{1}{c}{$<$ 1} \\
    Pigeon-07 & \multicolumn{1}{c}{1.5} & \multicolumn{1}{c}{$<$ 1} \\
    Pigeon-08 & \multicolumn{1}{c}{7.4} & \multicolumn{1}{c}{$<$ 1} \\
    Pigeon-09 & \multicolumn{1}{c}{88.6} & \multicolumn{1}{c}{$<$ 1} \\
    Pigeon-10 & \multicolumn{1}{c}{1381.4} & \multicolumn{1}{c}{$<$ 1} \\
    Pigeon-100 & \multicolumn{1}{c}{--} & \multicolumn{1}{c}{$<$ 1} \\
    Pigeon-1000 & \multicolumn{1}{c}{--} & \multicolumn{1}{c}{1.0} \\
    Pigeon-10000 & \multicolumn{1}{c}{--} & \multicolumn{1}{c}{1.8} \\
    Pigeon-100000 & \multicolumn{1}{c}{--} & \multicolumn{1}{c}{4.2} \\
    Pigeon-1000000 & \multicolumn{1}{c}{--} & \multicolumn{1}{c}{45.1} \\
    Pigeon-10000000 & \multicolumn{1}{c}{--} & \multicolumn{1}{c}{664.0} \\
    Pigeon-100000000 & \multicolumn{1}{c}{--} & \multicolumn{1}{c}{--} \\
    \bottomrule
    \end{tabular}%
  \label{tab:pigeontwo}%
\end{table}%

\begin{table}[h]
  \centering
  \caption{The performance of the space-indexed model on the Van-loading and Van-loading-X instances after one hour. 
  Each row represents 10 pseudo-randomly generated instances.}
    \begin{tabular}{lrrr}
    \addlinespace
    \toprule
    Dataset & Box Types  & Boxes  & Gap \\
            &            & (Mean) & (Mean) \\
    \midrule
    Van-loading-01 & 3     & 78    & 0.10\% \\
    Van-loading-02 & 5     & 92    & 0.31\% \\
    Van-loading-03 & 8     & 81    & 0.42\% \\
    Van-loading-04 & 10    & 73    & 0.32\% \\
    Van-loading-05 & 12    & 73    & 0.09\% \\
    Van-loading-06 & 15    & 70    & 0.42\% \\
    Van-loading-07 & 20    & 68    & 0.50\% \\
    Van-loading-08 & 30    & 71    & 0.99\% \\
    Van-loading-09 & 40    & 69    & 1.26\% \\
    Van-loading-10 & 50    & 73    & 1.10\% \\
    Van-loading-X-01 & 3     & 133   & 0.16\% \\
    Van-loading-X-02 & 5     & 155   & 0.16\% \\
    Van-loading-X-03 & 8     & 132   & 0.56\% \\
    Van-loading-X-04 & 10    & 119   & 0.87\% \\
    Van-loading-X-05 & 12    & 123   & 0.29\% \\
    Van-loading-X-06 & 15    & 115   & 0.36\% \\
    Van-loading-X-07 & 20    & 116   & 0.30\% \\
    Van-loading-X-08 & 30    & 122   & 0.29\% \\
    Van-loading-X-09 & 40    & 119   & 0.18\% \\
    Van-loading-X-10 & 50    & 121   & 0.28\% \\
    \bottomrule
    \end{tabular}
  \label{tab:vanloading}
\end{table}

\section{Conclusions and Future Work}
\label{sec:conclusions}

The space-indexed linear programming relaxation is strong, but employs a pseudo-polynomial numbers of variables,
  which presents a challenge to off-the-shelf solvers.
Future work could focus on the development of solvers specific to such relaxations. 
One could either implement a matrix-free interior point method exploiting the structure of the relaxation, or employ a column generation schema. 
The Dantzig-Wolfe decomposition could, perhaps, divide the problem into sub-problems pertaining to each type of box (\ref{no-overlap})
and a master problem enforcing the set-packing constraint (\ref{muconstraint1}).
Another interesting direction for future research is automated reformulations. Just as modern integer programming solvers often fail on small instances of the Chen/Padberg formulation, they fail on certain formulations of scheduling problems. It has been suggested \cite{Marecek2011} one could try to discretise such relaxations automatically, starting with the identification of what variables play the role of relative position in time or space, which 
could improve the performance of general-purpose integer linear programming solvers considerably. 
This may help to clarify the importance of discretisation in linear programming relaxations in general.

\paragraph{Acknowledgements}
This work has been first presented at the 7th ESICUP Meeting in Buenos Aires, Argentina, June 6-9, 2010.

\bibliographystyle{elsarticle-num}
\bibliography{mip-refs}

\appendix

\begin{algorithm}[tb]
\caption{
$f(D, n, L, x, y, z, t)$
}
\label{algo:genNonOverlap}
\begin{algorithmic}[1]
  \STATE \textbf{Input:} Discretisation as indices $D \subset \R^3$ of $\mu$ variables \eqref{eq:muvars}, number of boxes $n$, dimensions $L_{\alpha i} \in \R \quad \forall \alpha = x, y, z, i=1, \ldots, n$, 
indices $(x, y, z) \in D$, $t \in \{ 1, 2, \ldots, n\} $ of one scalar within $\mu$ \eqref{eq:muvars} 
  \STATE \textbf{Output:} Set of indices of $\mu$ variables \eqref{eq:muvars} to include in a set-packing inequality \eqref{no-overlap} \rule{0pt}{2.5ex}
  \vspace{1ex}
   \STATE Set $S \leftarrow \{ (x, y, z, t) \} $
   \FOR{ each other unit $(x', y', z') \in D, (x, y, z) \not = (x', y', z')$ } 
   \FOR{ each box type $t' \in \{1, 2, \ldots, n \} $ }  
    \IF{ box of type $t$ at $(x, y, z)$ overlaps box of type $t'$ at $(x', y', z')$, i.e., 
		$(x \le x' + L_{x't'} \le x + L_{xt}) \land (y \le y' + L_{y't'} \le y + L_{yt})$ \\
     \quad \quad $\land (z \le z' + L_{z't'} \le z + L_{zt}) $
} 
        \STATE $S \leftarrow S \cup \{ (x', y', z', t') \}$
    \ENDIF
   \ENDFOR
   \ENDFOR
   \RETURN S
\end{algorithmic}
\end{algorithm}

\end{document}